\def\mymode{i}  

\documentclass{amsart}

\usepackage{hyperref}
\usepackage{amsmath}
\usepackage{amssymb}
\usepackage{mathrsfs}
\usepackage{amsthm}
\usepackage[utf8]{inputenc}
\usepackage{svninfo}

\makeatletter

\ifx\Presentation\undefined
\fi

\def\newrefformat#1#2{%
  \@namedef{pr@#1}##1{#2}}
\def\fref#1{\@prettyref#1:}
\def\@prettyref#1:#2:{%
  \expandafter\ifx\csname pr@#1\endcsname\relax%
    \PackageWarning{prettyref}{Reference format #1\space undefined}%
    \ref{#1:#2}%
  \else%
    \csname pr@#1\endcsname{#1:#2}%
  \fi%
}

\newrefformat{eq}{\eqref{#1}}
\newrefformat{cha}{Chapter~\ref{#1}}
\newrefformat{sec}{Section~\ref{#1}}
\newrefformat{apx}{Appendix~\ref{#1}}
\newrefformat{tab}{Table~\ref{#1} on page~\pageref{#1}}
\newrefformat{fig}{Figure~\ref{#1} on page~\pageref{#1}}
\newrefformat{item}{\ref{#1}}

\ifx\thmnum\undefined
\newtheorem{thmnum}{}[section]
\fi

\newcommand{\mynewthm}[3][]{%
  \def\PARAM{#1}
  \ifx\PARAM\empty
  \newtheorem{#2}[thmnum]{#3}
  \else
  \newtheorem{#2}{#3}[#1]
  \fi
  \newtheorem*{#2*}{#3}%
  \newrefformat{#2}{#3~\ref{##1}}%
}

\ifx\undefined\undefined
\newcommand{\ThmLabel}{Theorem}
\newcommand{\PrpLabel}{Proposition}
\newcommand{\LemLabel}{Lemma}
\newcommand{\FctLabel}{Fact}
\newcommand{\CorLabel}{Corollary}
\newcommand{\DfnLabel}{Definition}
\newcommand{\ConvLabel}{Convention}
\newcommand{\NtnLabel}{Notation}
\newcommand{\CstLabel}{Construction}
\newcommand{\ExmLabel}{Example}
\newcommand{\RmkLabel}{Remark}
\newcommand{\QstLabel}{Question}

\else
\newcommand{\ThmLabel}{\iflanguage{french}{Théorème}{Theorem}}
\newcommand{\PrpLabel}{Proposition}
\newcommand{\LemLabel}{\iflanguage{french}{Lemme}{Lemma}}
\newcommand{\FctLabel}{\iflanguage{french}{Fait}{Fact}}
\newcommand{\CorLabel}{\iflanguage{french}{Corollaire}{Corollary}}
\newcommand{\DfnLabel}{\iflanguage{french}{Définition}{Definition}}
\newcommand{\ConvLabel}{Convention}
\newcommand{\NtnLabel}{Notation}
\newcommand{\CstLabel}{Construction}
\newcommand{\ExmLabel}{\iflanguage{french}{Exemple}{Example}}
\newcommand{\RmkLabel}{\iflanguage{french}{Remarque}{Remark}}
\newcommand{\QstLabel}{Question}

\fi

\theoremstyle{plain}
\mynewthm{thm}{\ThmLabel}
\mynewthm{prp}{\PrpLabel}
\mynewthm{lem}{\LemLabel}
\mynewthm{fct}{\FctLabel}
\mynewthm{cor}{\CorLabel}

\theoremstyle{definition}
\mynewthm{dfn}{\DfnLabel}
\mynewthm{conv}{\ConvLabel}
\mynewthm{conj}{Conjecture}
\mynewthm{ntn}{\NtnLabel}
\mynewthm{cst}{\CstLabel}

\theoremstyle{remark}
\mynewthm{rmk}{\RmkLabel}
\mynewthm{qst}{\QstLabel}
\mynewthm{exm}{\ExmLabel}

\ifx\Presentation\undefined

\newcommand{\myenumlabel}[1]{\textnormal{(\roman{#1})}}
\renewcommand{\theenumi}{\myenumlabel{enumi}}
\fi

\renewcommand{\today}{%
  \number\day\space
  \ifcase\month\or
  January\or February\or March\or April\or May\or June\or
  July\or August\or September\or October\or November\or December\fi
  \space \number\year}

\newcounter{cycprfcnt}
{\begin{list}{\PackageWarning{begnac}{Label required for cycprf}}%
  {%
    \setcounter{cycprfcnt}{1}
    \setlength{\itemindent}{0.5\leftmargin}%
    \setlength{\leftmargin}{0pt}%
    \newcommand{\cpcurr}{\myenumlabel{cycprfcnt}}%
    \newcommand{\cpnext}{\addtocounter{cycprfcnt}{1}\cpcurr}%
    \newcommand{\cpprev}{\addtocounter{cycprfcnt}{-1}\cpcurr}%
    \newcommand{\cpnum}[1]{\setcounter{cycprfcnt}{##1}\cpcurr}%
    \newcommand{\cpfirst}{\cpnum{1}}%
  }%
}%
{\qedhere\end{list}}%

\def\indsym#1#2{%
  \setbox0=\hbox{$\m@th#1x$}%
  \kern\wd0%
  \hbox to 0pt{\hss$\m@th#1\mid$\hbox to 0pt{$\m@th#1^{#2}$\hss}\hss}%
  \lower.9\ht0\hbox to 0pt{\hss$\m@th#1\smile$\hss}%
  \kern\wd0}

\def\nindsym#1#2{%
  \setbox0=\hbox{$\m@th#1x$}%
  \kern\wd0%
  \hbox to 0pt{\hss$\m@th#1\not$\kern1.4\wd0\hss}
  \hbox to 0pt{\hss$\m@th#1\mid$\hbox to 0pt{$\m@th#1^{#2}$\hss}\hss}%
  \lower.9\ht0\hbox to 0pt{\hss$\m@th#1\smile$\hss}%
  \kern\wd0}

\def\dotminussym#1#2{%
  \setbox0=\hbox{$\m@th#1-$}%
  \kern.5\wd0%
  \hbox to 0pt{\hss\hbox{$\m@th#1-$}\hss}%
  \raise.6\ht0\hbox to 0pt{\hss$\m@th#1.$\hss}%
  \kern.5\wd0}

\renewcommand{\emptyset}{\varnothing}
\renewcommand{\setminus}{\smallsetminus}

\newcommand{\rest}{{\restriction}}

\newcommand{\half}[1][1]{\hbox{$\frac{#1}{2}$}}

\newcommand{\cU}{\mathcal{U}}
\newcommand{\cV}{\mathcal{V}}

\newcommand{\bN}{\mathbb{N}}

\newcommand{\bR}{\mathbb{R}}

\makeatother

\DeclareMathOperator{\wap}{WAP}
\newcommand{\dyad}[2][n]{\hbox{$\frac{#2}{2^{#1}}$}}

\title{Reflexive representability and stable metrics}

\author{Itaï \textsc{Ben Yaacov}}
\address{Itaï \textsc{Ben Yaacov} \\
  Université Claude Bernard -- Lyon 1 \\
  Institut Camille Jordan \\
  43 boulevard du 11 novembre 1918 \\
  69622 Villeurbanne Cedex \\
  France}
\urladdr{\url{http://math.univ-lyon1.fr/~begnac/}}

\author{Alexander Berenstein}
\address{Alexander Berenstein \\
  Universidad de los Andes \\
  Departamento de Matemáticas \\
  Carrera 1 \# 18A-10, Bogotá \\
  Colombia ; and
  \newline\indent
  Université Claude Bernard -- Lyon 1 \\
  Institut Camille Jordan \\
  43 boulevard du 11 novembre 1918 \\
  69622 Villeurbanne Cedex \\
  France}
\email{\url{aberenst@uniandes.edu.co}}
\urladdr{\url{http://matematicas.uniandes.edu.co/~aberenst}}

\author{Stefano Ferri}
\address{Stefano Ferri
  Universidad de los Andes \\
  Departamento de Matemáticas \\
  Carrera 1 \# 18A-10, Bogotá \\
  Colombia}
\email{\url{stferri@uniandes.edu.co}}
\urladdr{\url{http://matematicas.uniandes.edu.co/~stferri/ferri.html}}

\thanks{The first author is supported by
  ANR chaire d'excellence junior THEMODMET (ANR-06-CEXC-007) and
  by Marie Curie research network ModNet.
  The second author is supported by the aforementioned ANR grant
  and by a research grant of the
  Facultad de Ciencias de la Universidad de los Andes.
  The third author thanks the
  Facultad de Ciencias de la Universidad de los Andes
  for the support received via the
  Proyecto Semilla ``Large groups and semigroups and their actions''.}

\if\mymode i
\renewcommand{\today}{%
  \number\day\space
  \ifcase\month\or
  January\or February\or March\or April\or May\or June\or
  July\or August\or September\or October\or November\or December\fi
  \space \number\year}
\svnInfo $Id: WAPStability.tex 927 2009-06-10 16:25:56Z begnac $
\thanks{\textit{Revision} {\svnInfoRevision} \textit{of} \today}
\else
\date{\today}
\fi



\begin{document}

\begin{abstract}

\noindent It is well-known that a topological group can be represented as a
group of isometries of a reflexive Banach space if and only if its topology is
induced by weakly almost periodic functions (see \cite{Shtern:CompactSemitopologicalSemigroups}, \cite{Megrelishvili:OperatorTopologies}
and \cite{Megrelishvili:TopologicalTransformations}). We show that for a metrisable group this is equivalent to
the property that its metric is uniformly equivalent to a stable metric in
the sense of Krivine and Maurey (see \cite{Krivine-Maurey:EspacesDeBanachStables}).
This result is used to give a partial negative answer to a problem
of Megrelishvili.
\end{abstract}

\maketitle

\section{Introduction}

In this paper we shall show that, for a topological group $G$ equipped
left-invariant metric, the topology
of $G$ is induced by weakly almost periodic functions if and only if its
metric is uniformly equivalent to a stable one. 
This result allows to give
simple proofs that several classical Banach spaces cannot be
represented as isometries of a reflexive Banach space.

We start by introducing some definitions and terminology.
Let $G$ be a topological group and let $C_b(G)$ denote
the algebra of continuous, bounded, 
complex-valued functions on $G$
equipped with the supremum norm. 

\begin{dfn}
  A function $f\in C_b(G)$ is said to be \emph{weakly almost periodic}
  if the set of left translates
  $\{f_x\colon x\in G\}$ (where $f_x(y)=f(xy)$)
  is a relatively weakly compact subset of $C_b(G)$.
\end{dfn}
The set of all weakly almost periodic functions on $G$ forms a
closed subalgebra of $C_b(G)$ which we shall denote by
$\wap(G)$. The functions in $\wap(G)$ can be characterised by the
following criterion due to Grothendieck
(see \cite{Grothendieck:CriteresDeCompacite} for a proof):

\begin{fct}[Grothendieck's Criterion]
  \label{fct:GrothCrit}
  A function
  $f\in C_b(G)$ is weakly almost periodic if and only if 
  for all sequences $(s_n)_n$, $(t_m)_m$ in $G$, we have
  \begin{gather*}
    \lim_{n \rightarrow \infty} \lim_{m \rightarrow \infty}
    f(s_nt_m)
    =
    \lim_{m \rightarrow \infty} \lim_{n \rightarrow \infty}
    f(s_nt_m)
  \end{gather*}
  whenever the limits exist. Equivalently, $f\in\wap(G)$ if
  and only if for all sequences 
  $(s_n)_n$, $(t_m)_m$ in $G$ and ultrafilters $\cU$, $\cV$ on 
  $\bN$ we have 
  \begin{gather*}
    \lim_{n,\cU} \lim_{m, \cV} f(s_nt_m)
    =
    \lim_{m,\cV} \lim_{n,\cU} f(s_nt_m).
  \end{gather*}
  We recall that, given a countable set $\{x_n:n\in\bN\}$ in
  a topological space and an ultrafilter $\cU$ on $\bN$, we write 
  $\lim_{n,\cU}x_n=L$ if for every neighbourhood $N$ of $L$ the
  set $\{n\in\bN:x_n\in N\}$ belongs to $\cU$.
\end{fct}

It is a classical result due to A.\ I.\ Shtern
(see \cite{Shtern:CompactSemitopologicalSemigroups})
that every compact semitopological group
(and so {\it a fortiori\/} every compact
topological group) can be represented as a group of
isometries of a reflexive Banach space
equipped with the weak operator topology.
We recall that the
weak operator topology on the space $B(E)$ of all continuous
linear endomorphisms of a Banach space $E$ is 
the weak topology induced by the maps
$\psi_{v,f}(T)=\langle f, Tv\rangle$; where $v \in E$, $f\in
E^\ast$,
and $\langle \cdot, \cdot \rangle$ is the dual pairing
between $E$ and $E^\ast$.

Another classical result states that
the weak operator topology on the
unitary group of a Hilbert space
coincides with its strong operator topology
(the strong topology on $B(E)$ is defined as the weak topology
induced by the maps
$\rho_v (T)=T(v)$, $v\in E$). 
The same result is true for the isometry group of any reflexive
Banach space (see \cite{Megrelishvili:OperatorTopologies}).

Moreover, on a compact semitopological semigroup every function on $C_b(G)$
is weakly almost periodic. 
The following theorem
(see \cite{Shtern:CompactSemitopologicalSemigroups}, \cite{Megrelishvili:OperatorTopologies} and \cite{Megrelishvili:ReflexivelyNotUnitarilyRepresentable}) 
relates the
richness of the algebra $\wap(G)$ to the property that 
$G$ can be embedded
into the isometry group of a reflexive Banach space. 

\begin{fct}\label{fct:ReflRepr}
  Let $G$ be a topological group.
  Then there is a topological isomorphism of $G$ into the
  isometry group of a reflexive Banach space with the weak
  (or the strong) operator topology if and only if 
  the topology of $G$ is induced by $\wap(G)$, that is if and only if $\wap(G)$
  separates the identity $e$ of
  $G$ from every closed subset not containing $e$.
\end{fct}

If a group satisfies the two equivalent properties of the previous theorem
it is said to be \emph{reflexively representable}.

When $G$ is locally compact, the left regular representation
$\lambda(g)(f(x))=f(g^{-1}x)$, $g \in G$, $f\in L^p(G,\mu)$, $p\ge
2$, where $\mu$ is the left Haar measure,
establishes a topological isomorphism of $G$ into the
isometry group
of the reflexive Banach space $L^p(G)$. Locally compact
groups are therefore reflexively representable (it can also be proved directly
that locally compact groups have enough weakly almost periodic functions to
separate points from closed sets). On 
the opposite extreme we
find the group $H_+[0,1]$ of all orientation preserving
homeomorphisms of the interval
$[0,1]$ with the topology of uniform convergence on compact
subsets. It was proved by 
Megrelishvili (\cite{Megrelishvili:EverySemitopologicalSemigroup})
that this group has no nonconstant weakly almost periodic function.

\begin{fct}
  \label{fct:WAPUnifCont}
  Let $G$ be a topological group, $f \in \wap(G)$.
  Then $f$ is both right and left uniformly continuous, meaning that
  for every $\varepsilon > 0$ there exists a neighbourhood
  $e \in U$ such that
  $|f(x)- f(xg)|,|f(x) - f(gx)| < \varepsilon$
  for every $x \in G$, $g \in U$.
\end{fct}
\begin{proof}
  \cite{Uspenskij:Compactification}, discussion following Theorem~5.3.
\end{proof}

We pass now to the definition of a stable metric.

\begin{dfn}
  We say that a metric $d(x,y)$ is \emph{stable} if for all 
  \emph{bounded sequences} $(s_n)_n$, $(t_m)_m$ in $G$
  (meaning there is an $M$ such that
  $d(e,s_n)\leq M$, $d(e,t_m)\leq M$ for all $n$ and all $m$),
  and for all ultrafilters $\cU$, $\cV$ on $\bN$ we have 
  \begin{gather*}
    \lim_{n,\cU} \lim_{m, \cV} d(s_n,t_m)
    =
    \lim_{m,\cV} \lim_{n,\cU} d(s_n,t_m).
  \end{gather*}
\end{dfn}

Thus, for a bounded semigroup with an invariant metric $d$,
the metric $d$ is stable if and only if the function
$f(x) = d(e,x)$ is weakly almost periodic.

In what follows we shall study the relation between stable metrics and
reflexive representability for metrisable groups.
In particular, we shall prove that for metrisable groups, being
reflexively representable is the same as having their metric
uniformly equivalent to a stable one.

The strategy of the proof will be to construct a stable pre-metric
(defined below)
using a Urysohn technique and then to compose it with a suitable
function to make it a stable metric.
In the next section we shall concentrate 
only on how to construct this function.

\section{Correction of triangle deficiency}

In this section we isolate a tool which has already been used by the
first author in
\cite{BenYaacov:Morley,BenYaacov:DefinabilityOfGroups}.
This tool allows one to extract, by a continuous manipulation,
a metric from a pre-metric, namely, from a binary function which would
be a metric if not for the fact that it fails the triangle inequality
(modulo the additional assumption of local continuity).

Recall that a function $f$ from a topological space $X$ to $\bR$ is
\emph{upper semi-continuous} if for every $t$ the set
$\{x\in X\colon f(x)<t\}$ is open in $X$.
Equivalently, if whenever $f(x) < t$ there exists a neighbourhood $U$
of $x$ such that $f\rest_U < t$.

\begin{dfn}
  Let $g\colon (\bR^+)^2 \to \bR^+$ be a symmetric,
  weakly increasing (i.e., non decreasing) function,
  satisfying $g(0,v) \leq v$ for all $v$.
  \begin{enumerate}
  \item We say that $g$ is a \emph{strong TD function} if in addition
    it is upper semi-continuous.
  \item We say that $g$ is a \emph{TD function} if we only know
    that for every $v < t$ there is $\delta > 0$ such that
    $g(\delta,v+\delta) < t$.
  \end{enumerate}
\end{dfn}

Notice that a weakly increasing function $g\colon (\bR^+)^2 \to \bR^+$
is upper semi-continuous if and only if
whenever $g(u,v)  < t$ there is $\delta > 0$ such that
$g(u+\delta,v+\delta) < t$.
Thus strong TD implies TD.
Conversely, every TD function gives rise to a strong TD function
as follows.

\begin{lem}
  \label{lem:StrongTD}
  Assume $g$ is a TD function and let 
  \begin{gather*}
    \tilde g(u,v) = \inf \{ g(u',v')\colon u'>u,v'>v \}.
  \end{gather*}
  Then $\tilde g$ is a strong TD function and
  $g \leq \tilde g$.
\end{lem}
\begin{proof}
  Clearly $\tilde g$ is symmetric and weakly increasing.
  Moreover, if $\tilde g(u,v) < t$ there must exist by definition
  $\delta > 0$ such that $g(u+2\delta,v+2\delta) < t$, whence
  $\tilde g\big( u+\delta, v+\delta \big) < t$.
  Thus $g$ is upper semi-continuous.
  By assumption, for every $v < t$ there is
  $\delta > 0$ such that
  $g(\delta,v+\delta) < t \Longrightarrow \tilde g(0,v) < t$.
  Thus $\tilde g(0,v) \leq v$, as desired.
  The inequality $g \leq \tilde g$ follows directly from the assumption
  that $g$ is weakly increasing.
\end{proof}

\begin{dfn}
  A \emph{pre-metric} on a space $X$ is a function
  $h\colon X^2 \to \bR^+$ satisfying:
  \begin{enumerate}
  \item \emph{Reflexivity:} $h(x,y) = 0 \Longleftrightarrow x = y$.
  \item \emph{Symmetry:} $h(x,y) = h(y,x)$.
  \end{enumerate}
  We say that a pre-metric is \emph{locally continuous} if for every
  $\varepsilon > 0$ there is $\delta > 0$ such that:
  \begin{gather*}
    h(x,y) < \delta \Longrightarrow |h(x,z)-h(y,z)| < \varepsilon.
  \end{gather*}
  We say that two pre-metrics $h_1$ and $h_2$ are
  \emph{uniformly equivalent} if for
  every $\varepsilon > 0$ there is $\delta > 0$ such that:
  \begin{gather*}
    h_i(x,y) < \delta \Longrightarrow h_j(x,y) < \varepsilon \qquad i,j \in \{1,2\}.
  \end{gather*}
\end{dfn}

A metric is always a pre-metric, and two metrics are uniformly
equivalent as pre-metrics if and only if  each is uniformly continuous
with respect to the other.
A pre-metric $h$ is a metric if and only if it satisfies the triangle inequality.
To every pre-metric $h$ we may attach a function which measures its
\emph{triangle deficiency (TD)}:
\begin{gather*}
  {\rm TD}_h(u,v) = \sup \{ h(x,z)\colon x,y,z \in X,
  h(x,y) \leq u, h(y,z) \leq v \}.
\end{gather*}
Indeed, $h$ satisfies the triangle
inequality if and only if ${\rm TD}_h(u,v) \leq u+v$.

\begin{lem}
  \label{lem:PreMetricTD}
  Let $h$ be a locally continuous pre-metric on $X$.
  Then ${\rm TD}_h$ is a TD function.
\end{lem}
\begin{proof}
  Clearly ${\rm TD}_h$ is symmetric, weakly increasing and satisfies
  ${\rm TD}_h(0,v) \leq v$.
  It remains to show that for every
  $v < t$ there is $\delta>0$ such that
  ${\rm TD}_h(\delta,v+\delta)<t$.
  By local continuity there is $\delta > 0$ such that
  $h(x,y) < \delta \Longrightarrow |h(x,z)-h(y,z)| < \frac{t-v}{2}$
  for all $x,y,z \in X$.
  We may also assume that $\delta <  \frac{t-v}{2}$.

  Assume now that
  $h(x,y) \leq \delta$ and
  $h(y,z) \leq v+\delta < \frac{t+v}{2}$.
  By choice of $\delta$ we have
  $h(x,z) < h(y,z) + \frac{t-v}{2} < t$.
  Thus ${\rm TD}_h(\delta,v+\delta) < t$, as desired.
\end{proof}

\begin{dfn}
  A \emph{correction function} for
  $g\colon (\bR^+)^2 \to \bR^+$ is a continuous weakly
  increasing function
  $f\colon \bR^+ \to \bR^+$ satisfying:
  \begin{enumerate}
  \item $f(t) = 0 \Longleftrightarrow t = 0$.
  \item $f \circ g(u,v) \leq f(u) + f(v)$.
  \end{enumerate}
\end{dfn}

\begin{lem}
  \label{lem:PreMetricCorrection}
  Let $X$ be a space, $h$ a locally continuous pre-metric on $X$
  and $f$ a correction function for ${\rm TD}_h$.
  Then $d_1 = f \circ h$ is a metric on $X$, uniformly equivalent to
  $h$.
\end{lem}
\begin{proof}
  To see that $f \circ h$ is a metric all we need to check is the triangle
  inequality.
  Let $x,y,z \in X$ and let
  $u = h(x,y)$, $v = h(y,z)$, $w = h(x,z)$.
  Then $w \leq {\rm TD}_h(u,v)$, whereby
  $f(w) \leq f \circ {\rm TD}_h(u,v) \leq f(u) + f(v)$, as desired.

  Now let $\varepsilon > 0$.
  The correction function $f$ is continuous and thus uniformly
  continuous on the compact $[0,1]$, so there exists $\delta > 0$ such that
  $t < \delta  \Longrightarrow f(t) < \varepsilon$.
  Conversely, $f(\varepsilon) > 0$ and
  $f(t) < f(\varepsilon) \Longrightarrow t < \varepsilon$.
  Thus $h$ and $f \circ h$ are uniformly equivalent.
\end{proof}

\begin{lem}
  \label{lem:ExistCorrection}
  Let $g$ be a TD function.
  Then there exists a correction function $f$ for $g$.
\end{lem}
\begin{proof}
  Let $g$ be a TD function and let $\tilde g$ be as in
  \fref{lem:StrongTD}.
  Then $\tilde g$ is strong TD and $g \leq \tilde g$.
  If we $\tilde g$ admits a correction
  function $f$ then $f \circ g \leq f \circ \tilde g$, whereby $f$ is a
  correction function for $g$ as well.
  It is therefore enough to prove that every strong TD function $g$
  admits a correction function.

  Let $D = \left\{ \dyad{k}\colon n \in \bN, 0 < k \leq 2^n \right\}$
  be the set of dyadic numbers in ${]}0,1]$.
  We shall proceed with a Urysohn style construction of $f$, choosing
  open sets $U_q \subseteq \bR^+$ for $q \in D$ satisfying
  $U_q \subseteq \{t\colon f(t) < q\}
  \subseteq \{t\colon f(t) \leq q\} \subseteq U_{q'}$
  where $q < q'$.
  Since we wish $f$ to be weakly increasing, each such open set will
  be of the form $[0,r_q{[}$.
  It is therefore enough to
  choose the right end-points $\{r_q\}_{q \in D}$.
  We proceed as in
  \cite[Lemma~2.19]{BenYaacov:Morley}, so that 
  for all $q,q' \in D$ we have:
  \begin{itemize}
  \item[(a)] $0 < r_q \leq q$,
  \item[(b)] If $q < q'$ then $r_q < r_{q'}$ and
    $g(r_q,r_{q'-q}) < r_{q'}$.
  \end{itemize}
  We choose $r_q$ for $q = \dyad{k} \in D$ by induction on $n$.
  For $n = 0$ we choose $r_1 = 1$, noticing that all the
  requirements hold.
  Assume now that $r_q$ has already been chosen for all
  $q \in D_n = \{\dyad{k}\colon 0 < k \leq 2^n\}$, and we wish to choose
  $r_q$ for $q \in D_{n+1} \smallsetminus D_n$, i.e., for
  $q = \dyad[n+1]{k}$, $0 < k < 2^{n+1}$ odd.
  It will be convenient to write
  $q^- = q - \dyad[n+1]{1}$ and
  $q^+ = q + \dyad[n+1]{1}$, noticing that
  $q^+ \in D_n$ and $q^- \in D_n\cup\{0\}$.

  In case $q \geq \dyad[n+1]{3}$ we have $q^\pm \in D_n$.
  Let:
  \begin{align*}
    \label{eq:one} \tag{$1$}
    & s_{q,q'} = \sup \{ s\leq 1\colon g(s,r_{q'}) < r_{q^- + q'} \},
    && q' \in D_n \cap [0,1-q^-] \\
    & s_q = \min\{ q, r_{q^+}, s_{q,q'}\colon q' \in D_n\cap[0,1-q^-]\}, \\
    & r_q = \frac{r_{q^-}+s_q}{2}.
  \end{align*}
  For $q' \in D_n \cap [0,1-q^-]$ we have
  $g(r_{q^-},r_{q'}) < r_{q^- + q'}$ by the induction hypothesis.
  Since $g$ is upper
  semi-continuous we obtain $r_{q^-} < s_{q,q'}$, and thus
  $r_{q^-} < r_q < s_q \leq r_{q^+}$.

  Let us now consider the case $q = \dyad[n+1]{1}$.
  For $q' \in D_{n+1} \cap [2q,1]$ we have already chosen
  $r_{q'}$ and we may define:
  \begin{align*}
    \label{eq:four} \tag{$2$}
    & s_{q,0} = \sup \{ s\leq 1\colon g(s,\half r_{2q}) < r_{2q} \}, \\
    \label{eq:five} \tag{$3$}
    & s_{q,q'} = \sup \{ s\leq 1\colon g(s,r_{q'}) < r_{q + q'} \},
    && q' \in D_{n+1} \cap [2q,1-q] \\
    & s_q = \min\{
    q, r_{2q}, s_{q,0}, s_{q,q'}\colon q' \in D_{n+1}\cap[2q,1-q]
    \}, \\
    & r_q = \half s_q.
  \end{align*}
  By assumption on $g$
  we have $g(0,\half r_{2q}) \leq \half r_{2q} < r_{2q}$, whereby
  $s_{q,0} > 0$.
  Similarly, if $q' \in D_{n+1} \cap [2q,1-q]$ then
  $g(0,r_{q'}) \leq r_{q'} < r_{q+q'}$, so again $s_{q,q'} > 0$.
  Thus $0 < r_q < s_q \leq r_{q^+}$.

  Let us check that (a) and (b) hold.
  Indeed, $r_q < s_q \leq q$ and
  $r_{q^-} < r_q < r_{q^+}$ (where $r_0 = 0$),
  ensuring that $q \mapsto r_q$ is strictly positive and strictly increasing
  on $D_{n+1}$.
  We are left with checking that if $q,q' \in D_{n+1}$, $q < q'$,
  then $g(r_q,r_{q'-q}) < r_{q'}$.
  Possibly exchanging $q$ with $q'-q$ we may assume that
  $q \in D_{n+1} \smallsetminus D_n$.
  Les us consider cases:
  \begin{enumerate}
  \item $q \geq \dyad[n+1]{3}$, $q'-q \in D_n$.
    Since $r_q < s_{q,q'-q}$ we have by \fref{eq:one}:
    \begin{gather*}
      g(r_q,r_{q'-q}) < r_{q^- + (q'-q)} < r_{q'}.
    \end{gather*}
  \item $q \geq \dyad[n+1]{3}$, $q'-q \in D_{n+1}\smallsetminus D_n$.
    Since $r_q < s_{q,(q'-q)^+}$ we have by \fref{eq:one}:
    \begin{gather*}
      g(r_q,r_{q'-q}) \leq g(r_q,r_{(q'-q)^+})  < r_{q^- + (q'-q)^+} = r_{q'}.
    \end{gather*}
  \item $q = \dyad[n+1]{1}$, $q' > 2q$.
    Since $r_q < s_{q,q'-q}$ we have by \fref{eq:five}:
    \begin{gather*}
      g(r_q,r_{q'-q}) < r_{q + (q'-q)} = r_{q'}.
    \end{gather*}
  \item $q = \dyad[n+1]{1}$, $q' = 2q$.
    Since $r_q < s_{q,0}$ and $r_q \leq \half r_{2q}$ we have by \fref{eq:four}:
    \begin{gather*}
      g(r_q,r_q) \leq g(r_q,\half r_{2q}) < r_{2q}.
    \end{gather*}
  \end{enumerate}
  We have thus verified that all the requirement hold and the
  inductive construction may proceed.

  We define $f\colon \bR^+ \to [0,1]$ by
  \begin{gather*}
    f(t) = \sup \{q \in D\colon r_q < t\}
    \qquad (\text{where } \sup \emptyset = 0),
  \end{gather*}
  so in particular, $f$ is continuous on the left.
  Since $q \mapsto r_q$ is strictly increasing, this definition of $f$
  agrees with
  \begin{gather*}
    f(t) = \inf \{q \in D\colon r_q > t\}
    \qquad (\text{where } \inf \emptyset = 1),
  \end{gather*}
  whereby $f$ is also continuous on the right, and thus continuous.
  It is clear that $f$ is weakly increasing
  (it is constant on $[1,\infty)$)
  and that $f(0) = 0$.
  On the other hand, if $f(t) = 0$ then $t < r_q \leq q$ for all
  $q \in D$, so $t = 0$.
  Finally, assume that $f(u) + f(v) < t$.
  If $t > 1$ then $f\circ g(u,v) \leq 1 < t$ directly.
  Otherwise, we have $f(u) < t_1$, $f(v) < t_2 = t-t_1$.
  By definition of $f$
  there are $q_1 < t_1$ and $q_2 < t_2$ in $D$ such that
  $r_{q_1} > u$, $r_{q_2} > v$.
  Then:
  \begin{gather*}
    g(u,v) \leq g(r_{q_1},r_{q_2}) < r_{q_1+q_2}
    \quad \Longrightarrow \quad
    f\circ g(u,v) \leq q_1 + q_2 < t.
  \end{gather*}
  Thus $f\circ g(u,v) \leq f(u) + f(v)$,
  and $f$ is a correction function for
  $g$ as desired.
\end{proof}

\begin{thm}
  \label{thm:TDCorrection}
  Let $X$ be a set, $h$ a pre-metric on $X$.
  Then $h$ is locally continuous if and only if it is uniformly
  equivalent to a metric on $X$.
  Moreover, this metric can be taken to be of the form $f \circ h$ where
  $f\colon \bR^+ \to \bR^+$ is continuous.
\end{thm}
\begin{proof}
  It is easy to check that if $h$ is uniformly equivalent to a metric
  $d$ on $X$ then it is locally continuous.
  Conversely, assume $h$ is locally continuous.
  Then ${\rm TD}_h$ is a TD function by \fref{lem:PreMetricTD}.
  By \fref{lem:ExistCorrection} there exists a correction function
  $f$ for ${\rm TD}_h$.
  Finally, by \fref{lem:PreMetricCorrection}, $f \circ h$ is a metric
  uniformly equivalent to $h$.
\end{proof}

\section{Stability and $\wap$-functions}

We start the section with a lemma which will allow us to consider only
bounded metrics in our proofs.

\begin{lem}
  \label{lem:BoundedMetric}
  Let $G$ be a metric group with a stable metric $d$,
  then $\delta(x,y):=\frac{d(x,y)}{1+d(x,y)}$ is also stable. 
\end{lem}
\begin{proof}
  Consider now two
  sequences, call them $(a_n)_n$ and $(b_m)_m$, of elements of $G$
  and let $\cU$ and $V$ be two non-principal ultrafilters on $\bN$.
  We want to show that 
  $$
  \lim_{n, \cU}\lim_{m, \cV}\delta(a_n,b_m)
  =
  \lim_{m, \cV}\lim_{n, \cU}\delta(a_n,b_m).
  $$

  Given $N\in \bN$ define $U_N:=\{n\in\bN : d(e,a_n)<N\}$ 
  and $V_N:=\{m\in\bN : d(e,b_m)<N\}$. 
  Let also $\varepsilon>0$ and let $L = 1/\varepsilon$,
  so $x\geq L$ implies $x/(1+x)>1-\varepsilon$.
  We proceed by cases:

  \textbf{Case 1:} 
  There is $N\in \bN$ such that 
  $U_N\in \cU$ and $V_N\in \cV$.

  Since $d(x,y)$ is stable, we get that 
  $$
  \lim_{n, \cU}\lim_{m, \cV} d(a_n,b_m)
  =
  \lim_{m, \cV}\lim_{n, \cU} d(a_n,b_m).
  $$
  On the other hand, since $f(x)=x/(1+x)$ is uniformly 
  continuous, we get 
  $$
  \lim_{n, \cU}\lim_{m, \cV}\delta(a_n,b_m)
  = \lim_{n, \cU}\lim_{m, \cV} f(d(a_n,b_m))
  = \lim_{m, \cV}\lim_{n, \cU} f(d(a_n,b_m))
  = \lim_{m, \cV}\lim_{n, \cU}\delta(a_n,b_m).
  $$ 

  \textbf{Case 2:} There is $N\in \bN$ such that 
  $U_N\in \cU$ but that for all $M>0$ we have $V_M\notin \cV$.

  For $n\in \bN$ we have
  \begin{gather*}
    \{m\in \bN\colon d(a_n,b_m)\geq L\}
    \supseteq
    \{m\in \bN\colon d(b_m,e)\geq L +d(e,a_n)\}
    \in \cV,
    \intertext{whereby}
    \lim_{m, \cV}\delta(a_n,b_m)\in [1-\varepsilon,1]
    \intertext{and thus}
    \lim_{n, \cU}\lim_{m, \cV}\delta(a_n,b_m)\in [1-\varepsilon,1].
  \end{gather*}
  On the other hand,
  assume $m \notin V_{N+L}$.
  Then
  \begin{gather*}
    \{n\in \bN\colon d(a_n,b_m) > L\}
    \supseteq 
    \{n\in \bN\colon d(a_n,e) < N\}
    \in \cU,
    \intertext{so}
    \lim_{n, \cU}\delta(a_n,b_m)\in [1-\varepsilon,1],
    \intertext{and since $(V_{N+L})^c\in \cV$:}
    \lim_{m, \cV} \lim_{n, \cU}\delta(a_n,b_m)\in [1-\varepsilon,1].
  \end{gather*}
  Since $\varepsilon>0$ is arbitrary, we obtain
  $$
  \lim_{n, \cU}\lim_{m, \cV}\delta(a_n,b_m)
  =1
  = \lim_{m, \cV}\lim_{n, \cU}\delta(a_n,b_m).
  $$

  \textbf{Case 3:}
  For all $N\in \bN$ we have
  $U_N\notin \cU$ and $V_N\notin \cV$.

  Just like in the previous case, we show that
  \begin{gather*}
    \lim_{n, \cU}\lim_{m, \cV}\delta(a_n,b_m)
    = 1
    = \lim_{m, \cV}\lim_{n, \cU}\delta(a_n,b_m).
    \qedhere
  \end{gather*}
\end{proof}

We are now ready to prove the main result of our paper,
namely, that a
metrisable group is reflexive representable if and only if 
its metric is uniformly equivalent to a stable metric.
We shall prove separately the two implications. 

Chaatit \cite{Chaatit:Representation}
proved that the additive group of any stable
separable Banach space can be represented as a group
of isometries of a separable reflexive Banach space.
First we  show, applying results by Shtern (\fref{fct:ReflRepr}),
that any stable group is reflexively representable.

\begin{prp}
  Assume that $G$ is a group equipped with a left-invariant stable
  metric $d$.
  Then $G$ is reflexively representable.
\end{prp}
\begin{proof}
  By \fref{lem:BoundedMetric} we can assume without loss of
  generality that $d$ is bounded.
  By Grothendieck's Criterion (\fref{fct:GrothCrit}), and
  since $d(x,y)=d(e,x^{-1}y)$ is a stable distance,
  the function $f(x) = d(e,x)$ is weakly almost periodic.
  Clearly it separates the identity from every closed set not
  containing it.
  By \fref{fct:ReflRepr}, $G$ is reflexively representable.
\end{proof}

Conversely,

\begin{thm}
  \label{thm:ReflReprStable}
  Let $G$ be a group equipped with a left-invariant metric, and
  assume that $G$ is reflexively representable.
  Then there is a left-invariant stable metric $\delta$ on $G$ 
  that is equivalent to $d$.
\end{thm}
\begin{proof}
  As in the previous result, by \fref{lem:BoundedMetric} we may assume
  that $d$ is bounded, and using the results proved in the previous
  section it suffices to produce a left-invariant stable locally
  continuous pre-metric uniformly equivalent to $d$.

  We now construct by induction a sequence of symmetric functions
  $f_n\colon G\longrightarrow [0,1]$ in $\wap(G)$
  (here, $f$ symmetric means that $f(g) = f(g^{-1})$ for all
  $g \in G$)
  and radii $r_n > 0$ such that the following conditions hold:
  \begin{itemize}
  \item[(a)] $f_n\rest_{B(e,r_{n+1})} = 0$,
    where $B(e,r)$ denotes the open ball of radius $r$ around $e$.
  \item[(b)] $f_n\rest_{G \setminus B(e,r_n)} = 1$,
  \item[(c)] $r_n \leq 2^{-n}$.
  \end{itemize}
  We start with $r_0 = 1$.
  Given $r_n$ we choose $f_n$ and $r_{n+1}$ as follows.
  First, by \fref{fct:ReflRepr} and our hypothesis,
  the algebra $\wap(G)$ separates points from closed sets.
  Therefore there exists a function
  $f_n\colon G\longrightarrow [0,1]$ in $\wap(G)$, such that
  $f_n[G\smallsetminus B(e,r_n)]=1$ and $f_n(e)=0$.
  Possibly replacing $f_n$ with
  $\min \bigl( f_n(x),f_n(x^{-1}) \bigr)$
  we may assume that $f_n$ is symmetric.
  Possibly replacing $f_n$ with $(2f_n - 1)^+$ we may further
  assume that $f_n$ is zero on a neighbourhood of $e$, and we may
  choose $r_{n+1}$ such that $f_n\rest_{B(e,r_{n+1})} = 0$.
  We may further require that $r_{n+1} \leq 2^{-n-1}$.

  Once we are done we define
  \begin{gather*}
    h(x,y):=\sum_{n=0}^\infty\frac{f_n(x^{-1}y)}{2^{n+1}}.
  \end{gather*}
  We claim that $h$ is a left-invariant stable locally continuous
  pre-metric uniformly equivalent to $d$:
  \\
  -- Left invariance, reflexivity $h(x,x)=0$ and symmetry
  $h(x,y)=h(y,x)$ follow immediately from the construction.
  \\
  -- Stability follows from the fact that $h$ is a
  uniform limit of functions in $\wap(G)$ which is a norm
  closed sub-algebra of $C_b(G)$.
  \\
  -- By \fref{fct:WAPUnifCont} $h$ is
  uniformly continuous (with respect to $d$).
  Conversely assume that $h(x,y) > 2^{-n}$.
  This means that $f_m(x^{-1}y) > 0$ for some $m < n$,
  so $x^{-1}y \notin B(e,r_n)$, i.e.,
  $d(x.y) \geq r_n$.
  Thus $d$ and $h$ are uniformly equivalent, and it follows that
  $h$ is locally continuous.
\end{proof}

\section{Some applications}
In this last section we shall use the characterisation of
reflexive representability we just proved to establish that
the additive group of several classical Banach spaces is not
reflexively representable.

Megrelishvili in \cite{Megrelishvili:TopologicalTransformations},
\cite{Megrelishvili:ReflexivelyNotUnitarilyRepresentable} and
\cite{Megrelishvili:EverySemitopologicalSemigroup}
and Glasner and Megrelishvili in
\cite[Problem 6.11]{Glasner-Megrelishvili:NewAlgebrasOfFunctions}
explicitely asked whether it was possible to find an abelian
group which is not reflexively representable. The problem
was solved directly in \cite{Ferri-Galindo:WAPCompactification}
using a modification of Raynaud's
proof \cite{Raynaud:EspacesDeBanachSuperstables}
of non uniform embeddability of $c_0$ into $\ell_2$.
The same result can also be deduced from two
facts: the space $c_0$ cannot be uniformly embedded into a reflexive
Banach space (see \cite{Kalton:CoarseAndUniformEmbeddings})
and a reflexively representable metric group always embeds uniformly
into a reflexive Banach space
(see \cite[Proposition~3.4]{Megrelishvili:OperatorTopologies}).
A similar result also holds
for the quasi-reflexive James' space $J$: it is not uniformly 
embeddable into a reflexive space (see
\cite{Kalton:CoarseAndUniformEmbeddings}) 
and thus it is not reflexively representable. 

In the light of our results it can now be proved
very easily that $c_0$ is not reflexively representable. 
Notice that this result 
also answers Question~6.12 of
\cite{Megrelishvili:TopologicalTransformations}
since we show that there exist quotients of
reflexively representable groups which are not
reflexively representable. This because
$c_0$ is, as every separable Banach space, a quotient of $\ell_1$
which is reflexively representable (in fact, it is even
unitarily representable).

\begin{thm}
  The additive group of the Banach space $c_0$ is not
  reflexively representable.
\end{thm}

\begin{proof}
  Suppose that $c_0$ is reflexively representable.
  Then its metric would be uniformly equivalent to a stable one, which
  is impossible by
  \cite[Théorème~5.1]{Raynaud:EspacesDeBanachSuperstables}.
\end{proof}

We use now our result to prove that Tsirelson's space $T$ is
not reflexively representable. Since $T$ is reflexive,
this gives a partial negative answer
to \cite[Question 6.9]{Megrelishvili:TopologicalTransformations} 
(see also \cite{Megrelishvili:ReflexivelyNotUnitarilyRepresentable}),
as it shows that there are metrisable groups which are uniformly
embedded into a reflexive Banach space (in fact, additive groups
of reflexive Banach spaces) which are not reflexively representable.

\begin{thm}
  Tsirelson's space $T$ is not reflexively representable.
\end{thm}

\begin{proof}
  By \cite[Théorème~4.1]{Raynaud:EspacesDeBanachSuperstables},
  if a Banach space admits a translation invariant stable metric
  which is uniformly
  equivalent to its norm-induced metric, then it contains a
  copy of $\ell_p$ for some $p$. Hence 
  Tsirelson's space $T$ does not admit such a metric and,
  since the metric produced in \fref{thm:ReflReprStable}
  is clearly translation invariant, Tsirelson's space
  is not reflexively representable.
\end{proof}
In a similar way it follows that Tsirelson-like spaces such as
Schlumprecht's space $\mathcal{S}_f$
(see \cite[Corollary 13.31]{Benyamini-Lindenstrauss:GeoremtricNonlinearFunctionalAnalysis}) are 
also not reflexively representable.

\providecommand{\bysame}{\leavevmode\hbox to3em{\hrulefill}\thinspace}
\providecommand{\MR}{\relax\ifhmode\unskip\space\fi MR }
\providecommand{\MRhref}[2]{%
  \href{http://www.ams.org/mathscinet-getitem?mr=#1}{#2}
}
\providecommand{\href}[2]{#2}

\end{document}